\documentclass[
12pt]{amsart}
\usepackage{color}
\usepackage{hyperref}
\usepackage{amssymb,amsmath}
\usepackage{amsthm}
\usepackage{amsfonts}
\usepackage {graphicx} 

 \def\BbbE{\mathbb E}

\def\cB{{\mathcal B}}

\def\bbE{{\mathbb E}}
\def\BbbN{\mathbb N}
\def\I{\mathbb I}

\def\bP{{\mathbf P}}
\def\bE{{\mathbf E}}
\def\cC{{\mathcal C}} 
\def\cF{{\mathcal F}}
\def\cB{{\mathcal B}}
\def\cT{{\mathcal T}}
\def\cW{{\mathcal W}}
\def\cL{{\mathcal L}}

\def\frt{{\mathfrak t}}
\def\frF{{\mathfrak F}}
\def\frS{{\mathfrak S}}
\def\frN{{\mathfrak N}}

\def\vecPi{\overrightarrow{\Pi}}

\def\one{{\mathbb I}}

       \def\vecxi{\overrightarrow{\xi}}
       \def\vecx{\overrightarrow{x}}
       \def\vecX{\overrightarrow{X}}

\numberwithin{equation}{section}

\theoremstyle{plain} 
\newtheorem{thm}[equation]{Theorem}

\newtheorem{lem}[equation]{Lemma}

\newtheorem{ass}[equation]{Assumptions}

\theoremstyle{definition}
\newtheorem{defn}[equation]{Definition}

\theoremstyle{remark}

\begin{document}

  \title{On the Masami Yasuda stopping game
  }
  \author{Krzysztof J.~Szajowski}
  \def\rightmark{On the Masami Yasuda stopping game}
  \def\leftmark{K.~J.~Szajowski}
  \date{\today}

%
%
%
  %
%
\address{Institute of Mathematics and Computer Science\\
 Wroc\l{}aw University of Technology\\ 
 Wybrze\.ze Wyspia\'nskiego 27\\ 
 PL-50-370 Wroc\l{}aw, Poland}
  \email{Krzysztof.Szajowski@pwr.wroc.pl}
  %

\begin{abstract}
The sero-sum stopping game for the stochastic sequences has been formulated in late sixties of the twenty century by Dynkin~\cite{dyn69:game}. The formulation had the assumption about separability of decision moment of the players which simplified the construction of the solution. Further research by Neveu~\cite{nev75:discrete} extended the model by admitting more general behaviour of the players and their pay--offs. In new formulation there is the problem with existence of the equilibrium. The proper approach to solution of the problem without restriction of former models was developed by Yasuda~\cite{yas85:randomized}. The results was crucial in these research. The author made often reference to the Yasuda's~\cite{yas85:randomized} result in his works (see~\cite{sza93:double,sza94:ZOR,sza95:SIAM}) as well as see results of others stimulated by this paper. Withal, in this note another stopping game model, developed by Yasuda with coauthors (see e.g.~\cite{kuryasnak80:majority} and \cite{szayas95:voting}) is recalled. The application of the model to an analysis of system of detectors shows the power of the game theory methods.
  
In the last part of the paper I would like to express my personal relation to the Masami Yasuda game.


\end{abstract}

\maketitle

\renewcommand{\baselinestretch}{1.2}
\section{\label{intro} Introduction}
The mathematical modelling of economic and engineering systems in stochastic environment leads to various mathematical optimization and game theory problems. If the decision problem relies on choice of intervention moments one can formulate the model of such case as the optimal stopping problem. If it is allowed to react more than once the approach depends on the number of decision makers and their aims. If there is one decision maker and two reactions (or fix number of possible moment of actions) we have the optimal two stopping (multiple stopping) problem. When there are two decision makers with their prescribed aims we usually treat the problem as the stopping game. The related models bring very subtle mathematical questions concerning the correctness of the model, possibility of inference about rational strategies, their realization and the existence of solution in the formulated mathematical model. In this note I would like to focus our attention of two group of models not very precisely defined. The first one is related to the existence of solution under mild assumption on the processes defining the payoffs in the zero-sum stopping game related to problem introduced by Dynkin~\cite{dyn69:game} (see the section~\ref{DGIntro}). The second group of the problems, which I have applied recently to modelling the sensor networks, developed by Yasuda with co-authors (see e.g.~\cite{kuryasnak80:majority} and \cite{szayas95:voting}), is devoted to multivariate stopping problem when there are decision makers having some interactions between each others (see the section~\ref{votingGintro}). 

\subsection{\label{DGIntro}The randomize strategies in Dynkin´s game}
\emph{E. B. Dynkin}~\cite{dyn69:game} presented the following problem: two players observe a stochastic sequence $X_n$, $n=1,2,... $. Each of them chooses a stopping time, say $\lambda$ (resp. $\mu)$ be the stopping time chosen by the first (resp. the second) player. It is additionally assumed that the first player can stop at odd and the second at even moments. The pay-off is then: $R(\lambda,\mu)= \bE X_{\lambda\wedge \mu}$. 
The player 1 seeks to maximize the expected pay-off, and the player 2 seeks to minimize it, it means that the solution is a pair $(\lambda^\star,\mu^\star)$ such that 
\begin{equation}\label{equ1}
R(\lambda,\mu^\star)\leq (\lambda^\star,\mu^\star)\leq (\lambda^\star,\mu).
\end{equation}

\emph{J. Neveu}~\cite{nev75:discrete} modified this problem as follows: there are three random sequences $(X_n,\cF_n)$, $(Y_n,\cF_n)$ and $(W_n,\cF_n)$ with 
\begin{ass}\label{ass1}
\begin{equation}
X_n\le W_n\le Y_n \text{ for each $n\in\BbbN$,}
\end{equation}
\end{ass}
and the pay-off equals:
\begin{equation}\label{equ2}
R(\lambda,\mu)=\bE\{X_\lambda \I_{\{\lambda< \mu\}}+W_\lambda \I_{\{\lambda=\mu\}}+Y_\mu \I_{\{\mu<\lambda\}}\},
\end{equation}
where $\lambda,\mu\in\frS$ are stopping times with respect to $\cF_n$.This problem has solution which is presented in \cite{nev75:discrete}. 

When the assumption~\ref{ass1} is not fulfilled then, in general there are no equilibrium in the set of stopping times with respect of observed processes $(X_n,\cF_n)$, $(Y_n,\cF_n)$ and $(W_n,\cF_n)$. It is \emph{M. Yasuda} who shown in \cite{yas85:randomized} that the mixed extension of this game has equilibrium without the assumption~\ref{ass1}. The mixed extension in this case means that the set of strategies (stopping times) is extended to include randomized stopping time. First, a finite horizon problem is considered. Next, the existence of the value and the equilibrium point in the infinite horizon problem with a discount factor is proved under some natural assumption concerning the integrability of the considered processes.

In that time there were many mathematicians doing research in stopping game (see e.g. the papers by Zabczyk~\cite{zab84:stopping}, Stettner~\cite{ste84:closedness}, Ohtsubo~\cite{oht86:neveu}). The Yasuda's paper~\cite{yas85:randomized} stimulated further research of Rosenberg, Solan and Vieille~\cite{rossolvie01:randomized} and Laraki and Solan~\cite{larsol05:zero} in the mixed extension of the stopping game for the processes with continuous parameter.

\subsection{\label{votingGintro}The stopping processes by voting procedure}
Let us consider $p$ person stopping game related to the observation of a Markov chain. Let  $(X_n,\frF_n,{\bP}_x)$, $n=0,1,2,\ldots $, be a homogeneous Markov chain defined on a probability space $(\Omega ,\frF, {\bf P})$ with state space $(\Bbb{E},\cB)$. The players are able to observe the Markov chain sequentially. At each moment $n$ their knowledge is represented by $\frF_n$. Each player has his own utility function $f_i: \Bbb{E}\rightarrow \Re $, $i=1,2,\ldots ,p,$ and at each moment $n$ each player declares separately his willingness to stop the observation of the process. The effective ends of the process and realization of the payoffs appears when a suitable subset of players agree to it. The aim of each player is to maximize their expected payoffs. In fact, the problem will be formulated as a $p$ person non-cooperative game with the concept of Nash equilibrium \cite {nas51:noncoop} as the solution. On the other hand, one can say  that the considered multilateral stopping procedure is based on sequential voting (cf \cite{fis71:majority}, \cite{home94:voting}, \cite{tay95:voting} for monotone rule  concept and the mathematics of voting).

Such model has been considered in mine and Yasuda paper \cite{szayas95:voting} and Ferguson~\cite{fer05:MR2104375}. Both papers were continuation of Masami Yasuda and his co-workers, Kurano and Nakagami research published in \cite{kuryasnak80:majority}, \cite{yasnakkur82:monotone}, \cite{yas95:explicite}, \cite{yas97:Markovdecision}. They have investigated the multilateral version of the optimal stopping problem for independent, identically distributed $p$ dimensional random vectors $\overline{X_n}$. The gain function of the $i$-th player is $X_n^i$ ($i$-th coordinate of $\overline{X_n}$). In \cite{kuryasnak80:majority} the following class of strategies is used.
\begin{enumerate}
\item Each player can declare to stop at any stage. 

\item The majority level $r$ $(1\leq r\leq p)$ is chosen by the players at
the beginning 
of the game.

\item During the sequential observation process, if the number of players declaring to stop is greater than or equal to the level $r$, the process must be stopped.
\end{enumerate}
This class of strategies is generalized in \cite{yasnakkur82:monotone} to monotone rules. Roughly speaking, a monotone rule is a $p$ variate, non-decreasing logical function defined on $\{0,1\}^p$. In both papers the problem is formulated as a $p$ person, non-cooperative game with concept of Nash point as a solution. Paper \cite{kuryasnak80:majority} generalizes the unanimity case, i.e. $p=r$ solved by Sakaguchi \cite{sak73:bivariate}. The motivation for the model considered is the secretary problem (see \cite{gimo66:maxseq} for the formulation of the problem). A solution of some bivariate version of the secretary problem is given in \cite{kuryasnak80:majority}. Presman and Sonin \cite{preson75:game} treat this problem with another set of strategies. They considered the model in which each player's decision does not affect the stopping of the process but his reward only. Sakaguchi \cite{sak78:bivariate} and Kadane \cite{kad78:reversibility} have solved a
multilateral sequential decision problem in which decisions whether to stop are made by the players alternately, instead of simultaneous decision under a monotone rule.

The recent paper on the voting stopping problem are \cite{nakkuryas00:jump}. 

\section{\label{SensorNET}Sensors' network and stopping games}
In~\cite{sza11:multi} the construction of the mathematical model for a multivariate surveillance system is presented. It is assumed that there is net $\frN$ of $p$ nodes which register (observe) signals modeled by discrete time multivariate stochastic process. At each node the state is the signal at moment $n\in \BbbN$ which is at least one coordinate of the vector $\vecx_n\in\bbE\subset\Re^m$. The distribution of the signal at each node has two forms and depends on \emph{a pure} or \emph{a dirty} environment of the node. The state of the system change dynamically. We consider the discrete time observed signal as $m\geq p$ dimensional process defined on the fixed probability space $(\Omega,\cF,\bP)$. The observed at each node process is Markovian with two different transition probabilities (see \cite{sarsza11:transition} for details). In the signal the visual consequence of the transition distribution changes at moment $\theta_i$, $i\in\frN$ is a change of its character. To avoid false alarm the confirmation from other nodes is needed. The family of subsets (coalitions) of nodes are defined in such a way that the decision of all member of some coalition is equivalent with the claim of the net that the disorder appeared. It is not sure that the disorder has had place. The aim is to define the rules of nodes and a construction of the net decision based on individual nodes claims. Various approaches can be found in the recent research for description or modelling of such systems (see e.g. \cite{tarvee08:sensor}, \cite{ragvee10:Markov}). The problem is quite similar to a pattern recognition with multiple algorithm when the fusions of individual algorithms results are unified to a final decision. The proposed solution will be based on a simple game and the stopping game defined by a simple game on the observed signals. It gives a centralized, Bayesian version of the multivariate detection with a common fusion center that it has perfect information about observations and \emph{a priori} knowledge of the statistics about the possible distribution changes at each node. Each sensor (player) will declare to stop when it detects disorder at his region. Based on the simple game the sensors' decisions are aggregated to formulate the decision of the fusion center. The  sensors' strategies are constructed as an equilibrium strategy in a non-cooperative game with a logical function defined by a simple game (which aggregates their decision). 

This approach uses the general description of such multivariate stopping games presented in the section~\ref{votingGintro}. The voting aggregation rules are relieved by the simple game (see  Ferguson~\cite{fer05:MR2104375}) and the underlining processes form Markov sequences (see~\cite{szayas95:voting}).

The model of disorder detection at each sensor are presented in the next section. It allows to define the individual pay-off of the players (sensors). It is assumed that the sensors are distributed in homogeneous way in the guarded area and the intruders behaviour are well modelled  by symmetric random walk. By these assumptions in the section~\ref{disorderONsensor} the \emph{a priori} distribution of the disorder moment at each node can be chosen in such a way that it gives the best model of the structure of sensors and the behaviour of intruder . The section~\ref{coopnobcoop} introduces the aggregation method based on a simple game of the sensors. The section~\ref{noncooperative} contains derivation of the non-cooperative game and existence theorem for equilibrium strategy. The final decision based on the state of the sensors is given by the fusion center and it is described in the section~\ref{strategiesOFsensors}. The natural direction of further research is formulated also in the same section. A conclusion and resume of an algorithm for rational construction of the surveillance system is included in  the section~\ref{finalremarks}.   

The extension of non-cooperative games to the case when the communication between player is allowed leads to various solutions concepts. The voting stopping game is interesting approach also in this direction of research.

\section{\label{disorderONsensor}Detection of disorder at sensors}
Following the consideration of Section~\ref{intro}, let us suppose that the process $\{\vecX_n,n\in\BbbN\}$, $\BbbN=\{0,1,2,\ldots\}$, is observed sequentially in such a way that each sensor, \emph{e.g.} $r$ (gets its coordinates in the vector $\vecX_n$ at moment $n$). By assumption, it is a stochastic sequence that has the Markovian structure given random moment $\theta_r$, in such a way that the process after $\theta_r$ starts from state $\vecX_{n\; \theta_r-1} $. The objective is to detect these moments based on the observation of $\vecX_{n\; \cdot}$ at each sensor separately. There are some results on the discrete time case of such disorder detection which generalize the basic problem stated by Shiryaev~in~\cite{shi61:detection} (see e.g. Brodsky and Darkhovsky~\cite{brodar93:nonparametr}, Bojdecki~\cite{boj79:disorder}, Poor and Hadjiliadis~\cite{poohad09:quickest}, Yoshida~\cite{yos83:complicated}, Szajowski~\cite{sza92:detection}) in various directions.

Application of the model for the detection of traffic anomalies in networks has been discussed by Tartakovsky et al.~\cite{tarroz06:intrusions}. The version of the problem when the moment of disorder is detected with given precision will be used here (see~\cite{sarsza11:transition}).

\subsection{\label{sformProblem} Formulation of the problem}
The observable random variables $\{\vecX_n\}_{n \in \BbbN}$ are consistent with the filtration $\mathcal{F}_n$ (or $\cF_n = \sigma(\vecX_0,\vecX_1,\ldots,\vecX_n)$). The random vectors $\vecX_n$ take values in $(\bbE, \mathcal{B})$, where $\bbE\subset\Re^m$. On the same probability space there are defined unobservable (hence not measurable with respect to $\cF_n$) random variables $\{\theta_r\}_{r=1}^m$ which have the geometric distributions:
\begin{eqnarray}
\label{rozkladyTeta}
\bP(\theta_r = j) = p_r^{j-1}q_r, \mbox{ $q_r=1-p_r \in (0,1)$, $j=1,2,\ldots$.}
\end{eqnarray}

The sensor $r$ follows the process which is based on switching between two, time homogeneous and independent, Markov processes $\{X_{rn}^i\}_{n \in \BbbN}$, $i=0,1$, $r\in\frN$ with the state space $(\bbE, \mathcal{B})$, both independent of $\{\theta_r\}_{r=1}^m$. Moreover, it is assumed that the processes $\{X_{rn}^i\}_{n \in \BbbN}$ have transition densities with respect to the $\sigma$-finite measure $\mu$, i.e., for any $B\in\cB$ we have  
\begin{eqnarray}\label{TransProbab}
\bP_x^{i}(X_{r1}^{i}\in B)&=&\bP(X_{r1}^{i}\in B|X_{r0}^{i}=x)=\int_Bf_x^{ri}(y)\mu(dy).
\end{eqnarray}
The random processes $\{X_{rn}\}$, $\{X_{rn}^0\}$, $\{X_{rn}^1\}$ and the random variables $\theta_r$ are connected via the rule: conditionally on $\theta_r = k$
\begin{eqnarray*}
X_{rn}&=&\left\{\begin{array}{ll}
X_{rn}^0,&\mbox{ if $k>n$,}\\
X_{r\;n+1-k}^1,&\mbox{ if $k\leq n$,}
\end{array}
\right.
\end{eqnarray*}
where $\{X_{rn}^1\}$ is started from $X_{r\;k-1}^0$ (but is otherwise independent of $X_{r\;\cdot}^0$).

For any fixed $d_r \in \{0,1,2,\ldots\}$ we are looking for the stopping time $\tau_r^{*}\in \cT$ such that
\begin{equation}
\label{PojRozregCiagowMark-Problem}
  \bP_x( | \theta_r - \tau_r^{*} | \leq d_r ) = \sup_{\tau \in \frS^X} \bP_x( | \theta_r - \tau | \leq d_r )
\end{equation}
where $\frS^X$ denotes the set of all stopping times with respect to the filtration
$\{\mathcal{F}_n\}_{n \in \BbbN}$. The parameters $d_r$ determines the precision level of detection and it can be different for too early and too late detection. These payoff functions measure the chance of detection of intruder.

\subsection{Construction of the optimal detection strategy}
In \cite{sarsza11:transition} the construction of $\tau^{*}$ by transformation of the problem to the optimal stopping problem for the Markov process $\vecxi$ has been made, such that   $\vecxi_{rn}=(\underline{\vecX}_{r\;n-1-d_r,n},\Pi_n$), where  $\underline{\vecX}_{r\;n-1-d_r,n}=(\vecX_{r\;n-1-d_r},\ldots,\vecX_{r\;n})$ and $\Pi_{rn}$ is the posterior process:
\begin{eqnarray}
    \Pi_{r0} &=& 0,\nonumber\\
    \Pi_{rn} &=& \bP_x\left(\theta_r \leq n \mid \mathcal{F}_n\right),\; n = 1, 2, \ldots  \nonumber
\end{eqnarray}
which is designed as information about the distribution of the disorder instant $\theta_r$. In this equivalent the problem of the payoff function for sensor $r$ is $h_r(\vecx_{r\;d_r+2},\alpha)$.

\section{\label{coopnobcoop} The aggregated decision via the cooperative game}
There are various methods combining the decisions of several classifiers or sensors. 
Each ensemble member contributes to some degree to the decision at any point of the sequentially delivered states. The fusion algorithm takes into account all the decision outputs from each ensemble member and comes up with an ensemble decision. When classifier outputs are binary, the fusion algorithms include the majority voting \cite{lamkrz94:majority}, \cite{LamSue97:pattern}, na\"{\i}ve Bayes combination~\cite{dompaz93:bayesian}, behavior knowledge space~\cite{huasue95:multiple}, probability approximation~\cite{kankimkim97:optimal} and singular value decomposition~\cite{mer99:correspondence}.

The majority vote is the simplest. The extension of this method is a simple game.        

\subsection{A simple game}
Let us assume that there are many nodes absorbing information and make decision if the disorder has appeared or not. The final decision is made in the fusion center which aggregates information from all sensors. The nature of the system and their role is to detect intrusion in the system as soon as possible but without false alarm. 

The voting decision is made according to the rules of \emph{a simple game}. Let us recall that a coalition is a subset of the players.  Let ${\cC}=\{C:C\subset \frN\}$ denote the class of all coalitions.  
\begin{defn}(see \cite{Owe95:MR1355082}, \cite{fer05:MR2104375})
\emph{A simple game} is coalition game having the characteristic function, $\phi(\cdot):\cC\rightarrow\{0,1\}$.      
\end{defn}
Let us denote $\cW=\{C\subset\frN:\phi(C)=1\}$ and ${\cL}=\{C\subset\frN:\phi(C)=0\}$. The coalitions in $\cW$ are called the winning coalitions, and those from $\cL$ are called the losing coalitions.
\begin{ass}
By assumption the characteristic function satisfies the properties:
\begin{enumerate}
\item\label{equat1} $\frN\in\cW$;
\item\label{equat2} $\emptyset\in \cL$;
\item\label{equat3} (the monotonicity): $T\subset S\in \cL$ implies $T\in \cL$.
\end{enumerate}
\end{ass}

\subsection{The aggregated decision rule}
When the simple game is defined and the players can vote presence or absence, $x_i=1$ or $x_i=0$, $i\in\frN$, of the intruder then the aggregated decision is given by the logical function
\begin{equation}\label{aggregateFUNCTION}
\pi(x_1,x_2,\ldots,x_p)=\sum_{C\in\cW}\prod_{i\in C}x_i\prod_{i\notin C}(1-x_i).
\end{equation}
For the logical function $\pi $ we have (cf \cite{yasnakkur82:monotone})
\[
\pi (x^1,\ldots ,x^p)=
x^i\cdot \pi (x^1,\ldots ,\stackrel{i}{\breve{1}}%
,\ldots,x^p)
+\overline{x}^i\cdot \pi (x^1,\ldots ,\stackrel{i}{\breve{0}}%
,\ldots ,x^p).
\]

\section{\label{noncooperative}A non-cooperative stopping game}
Following the results of the author and Yasuda~\cite{szayas95:voting} the multilateral stopping of a Markov chain problem can be described in the terms of the notation used in the non-cooperative game theory (see \cite{nas51:noncoop}, \cite{dresh81:games}, \cite{moulin}, \cite{Owe95:MR1355082}). 
Let $(\vecX_n,\frF_n,{\bP}_x)$, $n=0,1,2,\ldots ,N$, be a homogeneous Markov chain with state space $(\BbbE,\cB)$. The horizon can be finite or infinite. The players are able to observe the Markov chain sequentially. Each player has their utility function $f_i: \BbbE\rightarrow \Re $, $i=1,2,\ldots ,p$, such that ${\bE}_x|f_i(\vecX_1)|<\infty $. If process is not stopped at moment $n$, then each player, based on $\frF_n,$ can declare independently their willingness to stop the observation of the process.

\begin{defn}
{\rm (see \cite{yasnakkur82:monotone})} An individual stopping strategy of the player $i$ (ISS) is the sequence of random variables $\{\sigma _n^i\}_{n=1}^N$, where $\sigma_n^i:\Omega \rightarrow \{0,1\}$, such that $\sigma _n^i$ is $\frF_n$-measurable.
\end{defn}

The interpretation of the strategy is following. If $\sigma _n^i=1$ then player $i$ declares that they would like to stop the process and  accept the realization of $X_n$. Denote $\sigma ^i=(\sigma _1^i,\sigma _2^i,\ldots ,\sigma _N^i)$ and let $\frS^i$ be the set of ISSs of player $i$, $i=1,2,\ldots ,p$. Define 
\[
\frS=\frS^1\times \frS^2\times \ldots\times \frS^p. 
\]
The element $\sigma =(\sigma ^1,\sigma ^2,\ldots ,\sigma ^p)^T\in \frS$ will be called the stopping strategy (SS). The stopping strategy $\sigma \in \frS$ is a random matrix. The rows of the matrix are the ISSs. The columns are the decisions of the players at successive moments. The factual stopping of the observation process, and the players realization of the payoffs is defined by the stopping strategy exploiting $p$-variate logical function. Let $\pi :\{0,1\}^p\rightarrow \{0,1\}$. In this stopping game model the stopping strategy is the list of declarations of the individual players. The aggregate function $\pi$ converts the declarations to an effective stopping time.

\begin{defn}
A stopping time $\frt_\pi (\sigma )$ generated by the SS $\sigma \in \frS$ and the aggregate function $\pi $ is defined by 
\[
\frt_\pi (\sigma )=\inf \{1\leq n\leq N:\pi (\sigma _n^1,\sigma _n^2,\ldots ,\sigma _n^p)=1\}
\]
$(\inf (\emptyset )=\infty )$. Since $\pi $ is fixed during the analysis we skip index $\pi $ and write $\frt(\sigma )=\frt_\pi (\sigma )$. 
\end{defn}

We have $\{\omega \in \Omega : \frt_\pi (\sigma )=n\} =\bigcap\nolimits_{k=1}^{n-1}\{\omega \in \Omega : \pi (\sigma _k^1,\sigma
_k^2,\ldots,\sigma _k^p)=0\}\cap \{\omega \in \Omega :\pi (\sigma_n^1,\sigma _n^2,\ldots,\sigma _n^p)=1\}\in \frF_n$, then the random  variable $\frt_\pi (\sigma )$ is stopping time with respect to
$\{\frF_n\}_{n=1}^N$. 
For any stopping time $\frt_\pi (\sigma )$ and $i\in \{1,2,\ldots ,p\}$, let
\[
f_i(X_{\frt_\pi (\sigma )})=\left\{
\begin{array}{ll}
f_i(X_n) & \mbox{if }\frt_\pi (\sigma )=n\mbox{,} \\
\limsup_{n\rightarrow \infty }f_i(X_n) & \mbox{if }\frt_\pi (\sigma )=\infty 
\end{array}
\right.
\]
(cf \cite{shi}, \cite{szayas95:voting}). If players use SS $\sigma \in \frS$ and the individual preferences are converted to the effective stopping time by the aggregate rule $\pi $, then player $i$ gets $f_i(X_{\frt_\pi (\sigma )})$. 

Let ${}^{*}\!\sigma =({}^{*}\!\sigma ^1,{}^{*}\!\sigma ^2,\ldots ,{}^{*}\!\sigma ^p)^T$  be fixed SS. Denote 
\[
{}^{*}\!\sigma (i)=({}^{*}\!\sigma ^1,\ldots ,{}^{*}\!\sigma ^{i-1},\sigma ^i,{}^{*}\!\sigma ^{i+1},\ldots,{}^{*}\!\sigma^p)^T.
\]

\begin{defn}
\label{equdef}{\rm (cf. \cite{szayas95:voting})} Let the aggregate rule $\pi $ be fixed. The strategy
${}^{*}\!\sigma =({}^{*}\!\sigma ^1,{}^{*}\!\sigma ^2,\ldots ,{}^{*}\!\sigma ^p)^T\in \frS$ is an equilibrium strategy with respect to $\pi $ if for each $i\in \{1,2,\ldots ,p\}$ and any $\sigma ^i\in \frS^i$ we have
\begin{equation} 
{\bE}_xf_i(\vecX_{\frt_\pi ({}^{*}\!\sigma)})\geq  {\bE}_xf_i(\vecX_{\frt_\pi({}^{*}\!\sigma(i))}). 
\label{defequ}
\end{equation}
\end{defn}
The set of SS $\frS$, the vector of the utility functions $f=(f_1,f_2,\ldots, f_p)$ and the monotone rule $\pi $ define the non-cooperative game $\mathcal{G}$ = ($\frS$,$f$,$\pi$). 
The construction of the equilibrium strategy $ {}^{*}\!\sigma \in \frS$ in $\mathcal{G}$ is provided in \cite{szayas95:voting}. 
For completeness this construction will be recalled here. Let us define an individual stopping set on the state space. 
This set describes the ISS of the player. With each ISS of player $i$ the sequence of stopping events $D_n^i=\{\omega :\sigma _n^i=1\}$ combines. 
For each aggregate rule $\pi$ there exists the corresponding set value function $\Pi :\frF\rightarrow \frF$ such that $\pi (\sigma_n^1,\sigma _n^2,\ldots ,\sigma _n^p)= \pi \{\one_{D_n^1}, \one_{D_n^2},\ldots,\one_{D_n^p}\}= \one_{\Pi(D_n^1,D_n^2,\ldots,D_n^p)}$. 
For solution of the considered game the important class of ISS and the stopping events can be defined by subsets ${\it{C}}^i \in \mathcal{B}$ of the state space $\BbbE$. 
A given set ${\it{C}}^i\in\mathcal{B}$ will be called the stopping set for player $i$ at moment $n$ if $D_n^i= \{\omega :X_n\in {\it{C}}^i\}$ is the stopping event. 

For the logical function $\pi $ we have  
\[
\pi (x^1,\ldots ,x^p)=
x^i\cdot \pi (x^1,\ldots ,\stackrel{i}{\breve{1}}%
,\ldots,x^p)
+\overline{x}^i\cdot \pi (x^1,\ldots ,\stackrel{i}{\breve{0}}%
,\ldots ,x^p).
\]
It implies that for $D^i\in \frF$
\begin{equation}
\begin{array}{ll}
\Pi (D^1,\ldots ,D^p)=
& \{D^i\cap \Pi (D^1,\ldots,\stackrel{i}{\breve{\Omega}} , \ldots ,D^p)\} \\
& \cup \{\overline{D}^i\cap \Pi (D^1,\ldots ,\stackrel{i}{\breve{
\emptyset}},\ldots ,D^p)\}.
\end{array}
\label{decomposition}
\end{equation}

Let $f_i$, $g_i$ be the real valued, integrable (i.e. ${\bf E}_x|f_i(X_1)|<\infty $) function defined on $\BbbE$. For fixed $D_n^j$, $ j=1,2,\ldots ,p$, $j\neq i$, and ${\it{C}}^i\in \mathcal{B}$ define
\[
\psi ({\it{C}}^i)={\bf E}_x\left[f_i(X_1)\one_{{}^i\!D_1(D_1^i)}+ 
g_i(X_1)\one_{\overline{{}^i\!D_1(D_1^i)}}\right] 
\]
where
${}^i\!D_1(A)=\Pi (D_1^1,\ldots ,D_1^{i-1},A,D_1^{i+1},\ldots ,D_1^p)$ and 
$D_1^i=\{\omega :X_n\in {\it{C}}^i\}$. Let $a^{+}=\max \{0,a\}$ and $a^{-}=\min \{0,-a\}$. 
\begin{lem} $\label{optimal}$
Let $f_i$, $g_i$, be integrable and let ${\it{C}}^j\in \mathcal{B}$,
$j=1,2,\ldots,p$, 
$j\neq i$, be fixed. Then the set ${}^{*}\!{\it{C}}^i=\{x\in 
\BbbE:f_i(x)-g_i(x)\geq 0\}\in
\mathcal{B}$ is such that
\[
\psi ({}^{*}\!{\it{C}}^i)=\sup\limits_{{\it{C}}^i\it{\in }\mathcal{B} }\psi
({\it{C}}^i)
\]
and
\begin{eqnarray}
\psi ({}^{*}\!{\it{C}}^i)
& = & {\bf E}_x(f_i(X_1)-g_i(X_1))^{+}\one_{{}^{i}\!D_1(\Omega )}
\label{optset} \\
& & - {\bf E}_x(f_i(X_1)-g_i(X_1))^{-}\one_{{}^{i}\!D_1(\Omega )} +{\bf
E}_xg_i(X_1). 
\nonumber
\end{eqnarray}
\end{lem}

Based on Lemma \ref{optimal} we derive the recursive formulae defining the
equilibrium point and the equilibrium payoff for the finite horizon game. 

\subsection{The finite horizon game\label{finite}} 

Let horizon $N$ be finite. If the equilibrium strategy ${}^{*}\!\sigma $ exists, then we denote $v_{i,N}(x)={\bf E}_xf_i(X_{t({}^{*}\!\sigma )})$ the equilibrium payoff of $i$-th player when $X_0=x$. For the backward induction we introduce a useful notation. Let $\frS_n^i=\{\{\sigma_k^i\},k=n,\ldots ,N\}$ be the set of ISS for moments $n\leq k\leq N$ and $\frS_n=\frS_n^1\times \frS_n^2\times \ldots \times\frS_n^p$. The SS for moments not earlier than $n$ is ${}^n\!\sigma =({}^n\!\sigma ^1,{}^n\!\sigma ^2,\ldots ,{}^n\!\sigma^p) \in \frS_n$, where ${}^n\!\sigma ^i=(\sigma _n^i,\sigma _{n+1}^i,\ldots,\sigma_N^i)$.
Denote
\[
t_n=t_n(\sigma )=t(^n\sigma )=\inf \{n\leq k\leq N:\pi (\sigma _k^1,\sigma
_k^2,\ldots ,
\sigma _k^p)=1\}
\]
to be the stopping time not earlier than $n$. 

\begin{defn}
The stopping strategy ${}^{n*}\!\sigma =({}^{n*}\!\sigma ^1,{}^{n*}\!\sigma^2,\ldots,{}^{n*}\!\sigma ^p)$ is an equilibrium in $\frS_n$ if
\[
\begin{array}{ll}
 {\bf E}_x f_i(X_{t_n({}^{*}\!\sigma   )}) \geq 
 {\bf E}_x f_i(X_{t_n({}^{*}\!\sigma(i))}) & {\bf P}_x-\mbox{a.e.}
\end{array}
\]
for every $i\in \{1,2,\ldots ,p\}$, where \[
{}^{n*}\!\sigma (i)=({}^{n*}\!\sigma^1,\ldots,
{}^{n*}\!\sigma^{i-1},{}^n\sigma^i,
{}^{n*}\!\sigma^{i+1},\ldots,{}^{n*}\!\sigma ^p). \]
\end{defn}

Denote
\[
v_{i,N-n+1}(X_{n-1})={\bf E}_x[f_i(X_{t_n({}^{*}\!\sigma )})|\frF_{n-1}]
= {\bf E}_{X_{n-1}}f_i(X_{t_n({}^{*}\!\sigma )}). \]
At moment $n=N$ the players have to declare to stop and $v_{i,0}(x)=f_i(x)$. Let us assume that the process is not stopped up to moment $n,$ the players are using the equilibrium strategies ${}^{*}\!\sigma _k^i$, $i=1,2,\ldots ,p,$ at moments $k=n+1,\ldots ,N$. Choose player $i$ and assume that other players are using the equilibrium strategies ${}^{*}\!\sigma _n^j$, $j\neq i$, and player $i$ is using strategy $\sigma_n^i$ defined by stopping set ${\it{C}}^i$. Then the expected payoff $\varphi_{N-n}(X_{n-1},{\it{C}}^i) $ of player $i$ in the game starting at moment $n$, when the state of the Markov chain at moment $n-1$ is $X_{n-1\mbox{,}}$ is equal to
\[
 \varphi _{N-n}(X_{n-1},{\it{C}}^i)=
 {\bf E}_{X_{n-1}}\left[f_i(X_n)\one_{{}^{i*}\!D_n(D_n^i)}+ 
 v_{i,N-n}(X_n)
 \one_{\overline{{}^{i*}\!D_n(D_n^i)}}\right], 
\]
where
${}^{i*}\!D_n(A)=\Pi({}^{*}\!D_n^1,\ldots,{}^{*}\!D_n^{i-1},A, 
{}^{*}\!D_n^{i+1},\ldots,{}^{*}\!D_n^p)$. 

By Lemma \ref{optimal} the conditional expected gain $\varphi _{N-n}(X_{N-n}, 
{\it{C}}^i)$ attains the maximum on the stopping set
${}^{*}\!{\it{C}}_n^i=\{x\in
\bbE:f_i(x)-v_{i,N-n}(x)\geq 0\}$ and \setcounter{equation}{0}
\begin{equation}
\begin{array}{lll}
v_{i,N-n+1}(X_{n-1})
&=&{\bf E}_x[(f_i(X_n)-v_{i,N-n}(X_n))^{+} 
\one_{{}^{i*}\!D_n(\Omega)}|\frF_{n-1}] \\
& &- {\bf E}_x[(f_i(X_n)-v_{i,N-n}(X_n))^{-} 
\one_{{}^{i*}\!D_n(\emptyset)}|\frF_{n-1}] \\
& &+ {\bf E}_x[v_{i,N-n}(X_n)|\frF_{n-1}] 
\end{array} \label{valueatn}
\end{equation}
${\bf P}_x-$a.e..
It allows to formulate the following construction of the equilibrium strategy and the equilibrium value for the game $\mathcal{G}$.

\begin{thm}
In the game $\mathcal{G}$with finite horizon $N$ we have the following solution. 
\begin{description}
\item[(i)] The equilibrium value $v_i(x)$, $i=1,2,\ldots ,p$, of the game 
$\mathcal{G}$ can be calculated recursively as follows: \begin{enumerate}
\item $v_{i,0}(x)=f_i(x)$;
\item For $n=1,2,\ldots ,N$ we have ${\bf P}_x-$a.e. 
\small
\begin{eqnarray*}
v_{i,n}(x)={\bf E}_x[(f_i(X_{N-n+1})-v_{i,n-1}(X_{N-n+1}))^{+} 
\one_{{}^{i*}\!D_{N-n+1}(\Omega )}|\frF_{N-n}] \\ 
- {\bf E}_x[(f_i(X_{N-n+1})-v_{i,n-1}(X_{N-n+1}))^{-} 
\one_{{}^{i*}\!D_{N-n+1}(\emptyset )}|\frF_{N-n}] \\ 
+ {\bf E}_x[v_{i,n-1}(X_{N-n+1})|\frF_{N-n}], 
\end{eqnarray*}
for $i=1,2,\ldots ,p$.
\normalsize
\end{enumerate}
\item[(ii)] The equilibrium strategy ${}^{*}\!\sigma \in \frS$ is 
defined by the SS of the players ${}^{*}\!\sigma _n^i$, where
${}^{*}\!\sigma _n^i=1$ if
$ X_n\in {}^{*}\!{\it{C}}_n^i$, and
${}^{*}\!{\it{C}}_n^i=\{x\in \bbE : f_i(x)-v_{i,N-n}(x) \geq 0\}$,
$n=0,1,\ldots ,N$.
\end{description}

We have $v_i(x)=v_{i,N}(x)$, and ${\bf E}_xf_i(X_{t({}^{*}\!\sigma
)})=v_{i,N}(x)$, 
$i=1,2,\ldots ,p$.
\end{thm}

\section{Infinite horizon game\label{infinite}} 

In this class of games the equilibrium strategy is presented in Definition
\ref {equdef} but in class of SS 
\[
\frS_f^{*}=\{\sigma \in \frS^{*}:{\bf E}_xf_i^{-}(X_{t(\sigma )})<\infty
\quad
\mbox{ for every } \ x\in \bbE\mbox{, }i=1,2,\ldots ,p\}. 
\]
Let ${}^{*}\!\sigma \in \frS_f^{*}$  be an equilibrium strategy. Denote
\[
v_i(x)={\bf E}_xf_i(X_{t({}^{*}\!\sigma )}). 
\]

Let us assume that ${}^{(n+1)*}\!\sigma \in \frS_{f,n+1}^{*}$ is constructed and it is an equilibrium strategy. If players $j=1,2,\ldots ,p$, $j\neq i$, apply at moment $n$ the equilibrium strategies ${}^{*}\!\sigma _n^j$ , player $i$ the strategy $\sigma _n^i$ defined by stopping set ${\mathcal{C}}^{i}$ and ${}^{(n+1)*}\!\sigma $ at moments $n+1,n+2,\ldots$, then the expected payoff of the player $i$, when history of the process up to moment $n-1$ is known, is given by 
\[
  \varphi _n(X_{n-1},{\it{C}}^i)
  ={\bf E}_{X_{n-1}}\left[f_i(X_n) 
  \one_{{}^{i*}\!D_n(D_n^i)}+
  v_i(X_n)\one_{\overline{{}^{i*}\!D_n(D_n^i) }}\right], 
\]
where ${}^{i*}\!D_n(A)=\Pi ({}^{*}\!D_n^1,\ldots ,{}^{*}\!D_n^{i-1},A,{}^{*}\!D_n^{i+1},\ldots ,{}^{*}\!D_n^p)$, $ {}^{*}\!D_n^j=\{\omega \in \Omega:{}^{*}\!\sigma _n^j=1\}$, $j=1,2,\ldots ,p$,
$j\neq i$, and $D_n^i=\{\omega \in \Omega :\sigma _n^i=1\}=1\}= \{\omega \in \Omega :X_n\in \mathcal{C}^i\}$. By Lemma \ref{optimal} the conditional expected gain $\varphi _n(X_{n-1},{\it{C}}^i)$ attains the maximum on the stopping set ${}^{*}\!{\it{C}}_n^i=\{x\in \bbE:f_i(x)\geq v_i(x)\}$ and 
\begin{eqnarray*}
 \varphi _n(X_{n-1},{}^{*}\!{\it{C}}^i)
 &=&{\bf E}_x[(f_i(X_n)-v_i(X_n))^{+}
 \one_{{}^{i*}\!D_n(\Omega )}|\frF_{n-1}] \\
 &&- {\bf E}_x[(f_i(X_n)-v_i(X_n))^{-}
 \one_{{}^{i*}\!D_n(\emptyset )}|\frF_{n-1}] \\
 &&+ {\bf E}_x[v_i(X_n)|\frF_{n-1}].
\end{eqnarray*}

Let us assume that there exists solution $(w_1(x),w_2(x),\ldots ,w_p(x))$
of the equations
\setcounter{equation}{0}
\begin{eqnarray}
w_i(x) &=&
{\bf E}_x(f_i(X_1)-w_i(X_1))^{+}\one_{{}^{i*}\!D_1(\Omega )}
\label{equvalue} \\
&&-{\bf E}_x(f_i(X_1)-w_i(X_1))^{-}\one_{{}^{i*}\!D_1(\emptyset )} 
+ {\bf E}_xw_i(X_1),
\nonumber
\end{eqnarray}
$i=1,2,\ldots ,p$. Consider the stopping game with the following payoff
function for 
$i=1,2,\ldots ,p$.
\[
\phi _{i,N}(x)=\left\{
\begin{array}{ll}
f_i(x) & \mbox{if }n<N, \\
v_i(x) & \mbox{if }n\geq N.
\end{array}
\right.
\]

\begin{lem}\label{auxinfgame}
Let ${}^{*}\!\sigma {}\in \frS_f^{*}$ be an equilibrium strategy in the infinite horizon game $\mathcal{G}$. For every $N$ we have 
\[
{\bf E}_x\phi _{i,N}(X_{t^{*}})=v_i(x).
\]
\end{lem}

Let us assume that for $i=1,2,\ldots ,p$ and every $x\in \bbE$ we have 
\begin{equation}
{\bf E}_x[\sup\nolimits_{n\in \BbbN}f_i^{+}(X_n)]<\infty . \label{supplus} 
\end{equation}

\begin{thm}
Let $(X_n,\frF_n,{\bf P}_x)_{n=0}^\infty $ be a homogeneous Markov chain and 
the payoff functions of the players fulfill (\ref{supplus}). If $t^{*}=t({}^{*}\!\sigma )$, 
${}^{*}\!\sigma \in \frS_f^{*}$ then ${\bf E}_xf_i(X_{t^{*}})=v_i(x)$.
\end{thm}

\begin{thm}
Let the stopping strategy ${}^{*}\!\sigma \in \frS_f^{*}$ be defined by the 
stopping sets
${}^{*}\!{\it{C}}_n^i=\{x\in \bbE:f_i(x)\geq v_i(x)\}$, $i=1,2,\ldots ,p$, 
then ${}^{*}\!\sigma $ is the equilibrium strategy in the infinite stopping
game 
$\mathcal{G}$.
\end{thm}

\subsection{\label{strategiesOFsensors}Determining the strategies of sensors}
Based on the model constructed in Sections~\ref{disorderONsensor}--\ref{noncooperative} for the net of sensors with the fusion center determined by a simple game, one can determine the rational decisions of each nodes. The rationality of such a construction refers to the individual aspiration for the highest sensitivity to detect the disorder without false alarm. The Nash equilibrium fulfills requirement that nobody deviates from the equilibrium strategy because its probability of detection will be smaller. The role of the simple game is to define wining coalitions in such a way that the detection of intrusion to the guarded area is maximal and the probability of false alarm is minimal. The method of constructing the optimum winning coalitions family is not the subject of the research in this article. However, there are some natural methods of solving this problem.    

The research here is focused on constructing the solution of the non-cooperative stopping game as to determine the detection strategy of the sensors. To this end, the game analyzed in Section~\ref{noncooperative} with the payoff function of the players defined by the individual disorder problem formulated in Section~\ref{disorderONsensor} should be derived. 

The proposed model disregards correlation of the signals. It is also assumed that the fusion center has perfect information about signals and the information is available at each node. The further research should help to qualify these real needs of such models and to extend the model to more general cases. In some type of distribution of sensors, e.g. when the distribution of the pollution in the given direction is observed, the multiple disorder model should work better than the game approach. In this case the \emph{a priori} distribution of disorder moment has the form of sequentially dependent random moments and the fusion decision can be formulated as the threshold one: stop when $k^*$ disorder is detected. The method of a cooperative game was used in \cite{ghakri10:cooperativeMR2780188} to find the best coalition of sensors in the problem of the target localization. The approach which is proposed here shows possibility of modelling the detection problem by multiple agents at a general level. 

\subsection{\label{finalremarks}Final conclusion concerning the disorder detection system}
In a general case the consideration of the paper~\cite{sza11:multi} leads to the algorithm of constructing the disorder detection system. 
\begin{figure}[h]
\centerline{\includegraphics[width=2.8in,height=4in]{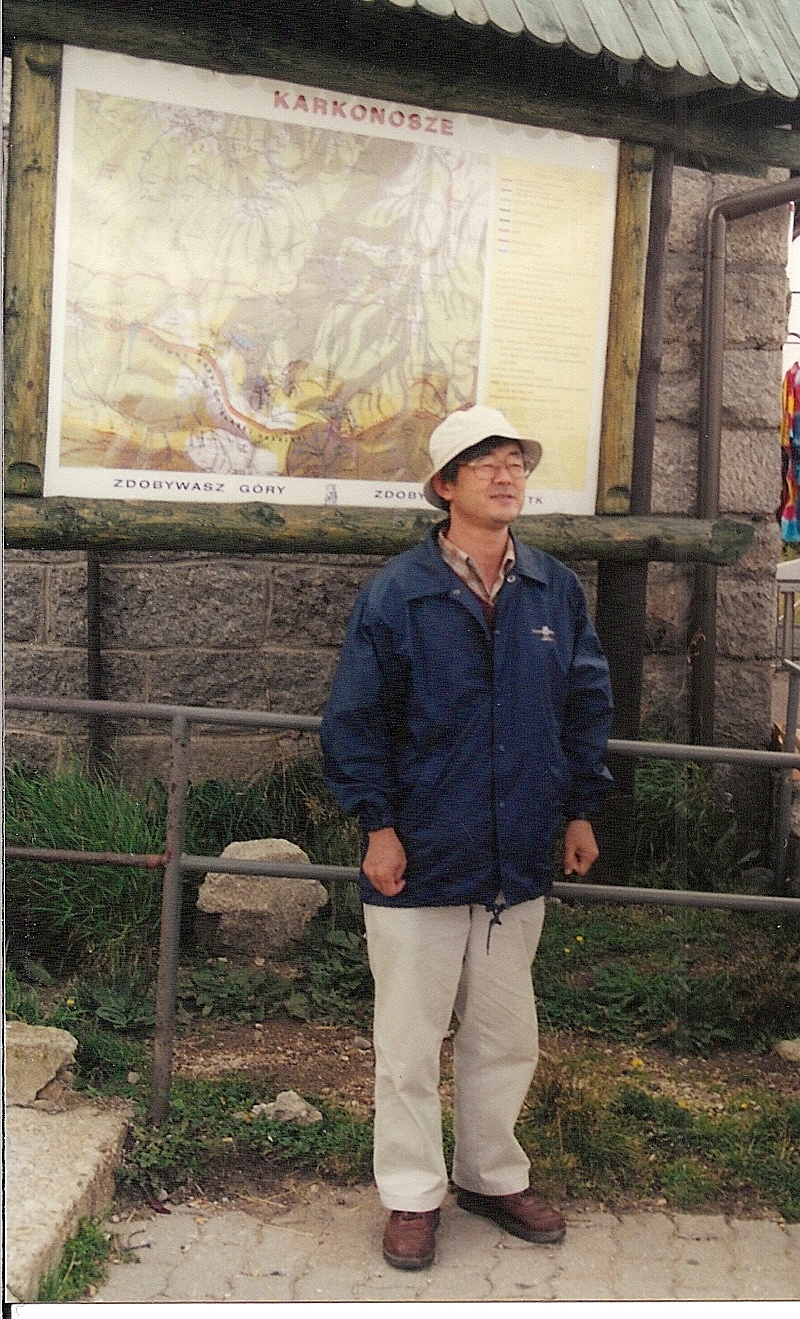}} 
\caption{\label{sudety2000}Sudety Mountains. International Conference on Mathematical Statistics STAT'2000. Szklarska Por\c{e}ba, Poland, August, 2000}
\end{figure}
\subsubsection{Algorithm}
\begin{enumerate}
\item Define a simple game on the sensors.
\item Describe signal processes and \emph{a priori} distribution of the disorder moments at all sensors. Establish the \emph{a posteriori} processes: $\vecPi_n=(\Pi_{1n},\ldots,\Pi_{mn})$, where $\Pi_{kn}=\bP(\theta\leq n|\cF_n)$. 
\item Solve the multivariate stopping game on the simple game to get the individual strategies of the sensors.
\end{enumerate}
\section{\label{final}The contribution to the mathematical education, the scientific cooperation and the friendship }
\begin{figure}[h]
\centerline{\includegraphics[width=3.2in,height=4in]{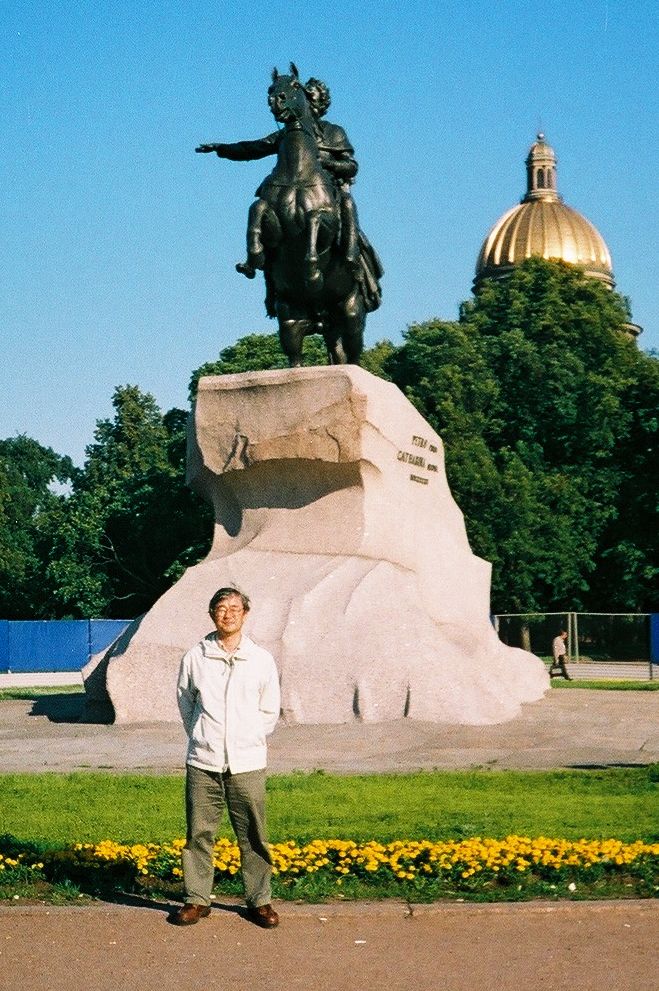}} 
\caption{\label{petersburg2008} 13${}^\text{th}$ Symposium of ISDG  2008, St.Petersburg, Russia}
\end{figure}

It was 1994 when I came to Japan for the first time based on Prof. Minoru Sakaguchi and Prof. Katsunori Ano invitation to take part in the International Conference on Stochastic Models and Optimal Stopping, Nanzan University, Nagoya. Since this event Professor Masami Yasuda is my guide in the mathematics and the Japanese culture. When the Internet connected the people we discussed the game model, which I call myself the Masami Yasuda game described in this note in  the sections~\ref{votingGintro} and \ref{SensorNET}. Based on the discussion we have written the papers~\cite{szayas95:voting} and \cite{yassza02:dynkin}.  
\begin{figure}
\centerline{\includegraphics[width=4in,height=2.8in]{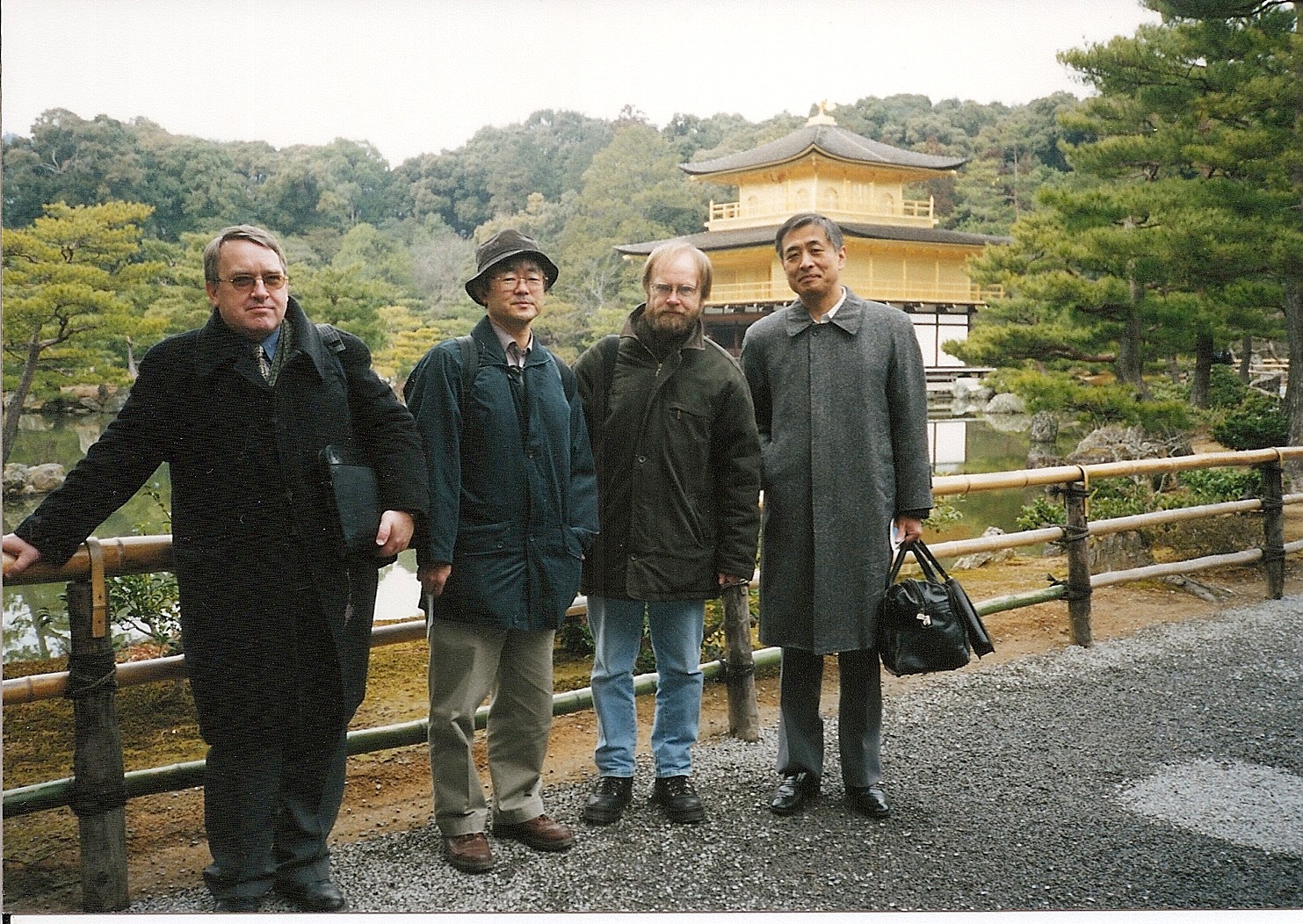}} 
\caption{\label{kyoto2002}RIMS Conference 2002, Kyoto, Japan. From the left: the author of this note, professors Masami Yasuda, Vladimir V. Mazalov and Mitsushi Tamaki}
\end{figure}
Our meeting and discussion were in Poland (see Figure~\ref{sudety2000}), Japan~(see Figure~\ref{kyoto2002}) and Russia~(see Figure~\ref{petersburg2008}). Last year, during the one day workshop, the possible further research and academic cooperation was the topic of our discussion.  
Based on European Union integration processes the Polish educational system is under very intensive reconstruction process. The mathematical education of engineering faculties students is very fragile task. Professor Yasuda provides us his extensive academic experience by his contribution to our local conferences devoted to teaching mathematics for non-mathematical major students. Such organized events were in Ibaraki National College of Technology, Hitachinaka (2006) and Wroc\l{}aw Univeristy of Technology (2008).

\section*{Acknowledgments}
The authors gratefully acknowledge the many helpful suggestion of the anonymous referees and the editor during the preparation of the paper.

\noindent\vspace{-18pt}
\def\cprime{$'$} \def\cprime{$'$} \def\cprime{$'$}
\providecommand{\bysame}{\leavevmode\hbox to3em{\hrulefill}\thinspace}
\providecommand{\MR}{\relax\ifhmode\unskip\space\fi MR }
\providecommand{\MRhref}[2]{%
  \href{http://www.ams.org/mathscinet-getitem?mr=#1}{#2}
}
\providecommand{\href}[2]{#2}

\end{document}